\documentclass[11pt]{amsart}

\usepackage{amsmath}
\usepackage{amssymb}
\usepackage[all]{xy}
\usepackage{sty1}

\makeatletter
\@addtoreset{equation}{section}
\makeatother

\def\bltwo#1{\mathrm{L}^2(#1)}
\def\falg{\mathrm{A}(G)}
\def\fal#1{\mathrm{A}(#1)}
\def\falig{\mathrm{A}_I(G)}
\def\fsal#1{\mathrm{B}(#1)}
\def\fsalg{\mathrm{B}(G)}
\def\fsalh{\mathrm{B}(H)}
\def\fsali#1{\mathrm{B}_I(#1)}
\def\fsalig{\mathrm{B}_I(G)}
\def\fsalih{\mathrm{B}_I(H)}
\def\faltg{\mathrm{A}_\tau(G)}
\def\falb#1#2{\mathrm{A}_{#1}(#2)}
\def\falbg#1{\mathrm{A}_{#1}(G)} 
 
\def\hdtdg{\mathfrak{HD}_{td}(G)}
\def\idemg{\mathrm{Idem}(G)}
\def\idemcg{\mathrm{Idem}_c(G)}
\def\spineg{\mathrm{A}^*(G)}
\def\tg{\mathcal{T}(G)}
\def\tnq#1{\mathcal{T}_{nq}(#1)}
\def\tnqg{\mathcal{T}_{nq}(G)}
\def\tnqtd#1{\mathcal{T}_{nq}^{td}(#1)}
\def\tnqtdg{\mathcal{T}_{nq}^{td}(G)}
\def\ttd#1{\mathcal{T}_{td}(#1)}
\def\ttdg{\mathcal{T}_{td}(G)}
\def\tnqg{\mathcal{T}_{nq}(G)}
\def\vnig{\mathrm{VN}_I(G)}
\def\vnih{\mathrm{VN}_I(H)}

\def\vntg{\mathrm{VN}_\tau(G)}

\def\genlin#1{\mathrm{GL}(#1)}
\def\Oh{\mathbb{O}}

\begin{document}

\newtheorem{falgidem}{Proposition}[section]

\newtheorem{ttdgisttdge}{Proposition}[section]
\newtheorem{tdcompact}[ttdgisttdge]{Proposition}
\newtheorem{tnqtdgchar}[ttdgisttdge]{Theorem}
\newtheorem{tnqtdgchar1}[ttdgisttdge]{Corollary}
\newtheorem{fsalgidem}[ttdgisttdge]{Theorem}
\newtheorem{spectrum}[ttdgisttdge]{Theorem}

\newtheorem{cbhomo}{Theorem}[section]
\newtheorem{isomorphism}[cbhomo]{Theorem}
\newtheorem{amenability}[cbhomo]{Proposition}

\title[The Algebra Generated by Idempotents]
{The Algebra Generated by Idempotents in a Fourier-Stieltjes Algebra}

\author{Monica Ilie and Nico Spronk}

\begin{abstract}
We study the closed algebra $\fsalig$ generated by the idempotents in the 
Fourier-Stieltjes algebra of a locally compact group $G$.  We show that it is
a regular Banach algebra with computable spectrum $G^I$, which we call
the idempotent compactification of $G$.  For any locally
compact groups $G$ and $H$, we show that $\fsalig$ is completely
isometrically isomorphic to $\fsalih$ exactly when $G/G_e\cong H/H_e$,
where $G_e$ and $H_e$ are the connected components of the identities.
We compute some examples to illustrate out results.
\end{abstract}

\maketitle

\footnote{{\it Date}: \today.

2000 {\it Mathematics Subject Classification.} Primary 43A30,43A25;
Secondary 43A22, 46L07, 43A70.
{\it Key words and phrases.} Fourier algebra, Fourier-Stieltjes algebra,
locally precompact group topology.

Monica Ilie's research supported by NSERC under grant no.\ 312585-05.

Nico Spronk's research supported by NSERC under grant no.\ 312515-05.}

\section{Introduction}

Let $G$ be a locally compact group and $\fsalg$ denote the {\it 
Fourier-Stieltjes algebra} of $G$, defined by Eymard \cite{eymard}.
Let $\idemg$ denote the family of all idempotents in $G$.  The theorem
of Host \cite{host} tells us exactly what each element $u$ of $\idemg$
is:  $u=1_Y$, the indicator function of an element $Y$ of $\Omega_o(G)$,
the Boolean ring of subsets generated by open cosets in $G$.

Since $\fsalg$ is commutative, $\idemg$ is a multiplicative semigroup 
in $\fsalg$, hence the algebra it generates is simply its linear span.  We let
$\fsalig$ denote the closure of the span of $\idemg$ in $\fsalg$.  This algebra
will be our main object of study.

Let us begin with a much simpler object.  Let $\falg$ denote the
{\it Fourier algebra} of $G$, also defined in \cite{eymard}.  
Then $\falg$ is an algebra of continuous 
functions on $G$ each of which vanish at infinity.  Moreover, evaluation 
functionals at points in $G$ comprise the Gelfand spectrum of $\falg$.  Let
$\falig$ denote the algebra generated by idempotent elements of $\falg$.

The following proposition is well-known.  
Our proof of (i) is very similar to the
one offered for abelian groups in \cite[8.1.4]{taylor}.

\begin{falgidem}\label{prop:falgidem}
{\bf (i)} Let $\idemcg$ denote the set of compactly
supported idempotents in $\fsalg$ and let $\Omega_o^c(G)$ denote the
Boolean ring of subsets generated by compact open cosets in $G$.  Then 
\[
\falg\cap\idemg=\idemcg=\{1_Y:Y\in\Omega_o^c(G)\}.
\]

{\bf (ii)} If $G_e$ denotes the connected component of the identity in $G$, 
then
$\falig=\fal{G\!:\!G_e}$, the subalgebra of all functions in $\falg$ which
are constant on cosets of $G_e$.
\end{falgidem}

\proof {\bf (i)}  If $u\in\falg\cap\idemg$ with $u\not=0$,
then there is an open compact
subset $Y$ of $G$ for which $u=1_Y$.  Moreover, any compactly supported
element of $\fsalg$ is in $\falg$, so $\falg\cap\idemg=\idemcg$.
It remains to show that $Y\in\Omega_o^c(G)$.
We note that for any $y\iin Y$, $y^{-1}Y$
is a compact open neighbourhood of $e$, so $G_e$ is compact.

Suppose first that $G_e=\{e\}$, so $G$ is totally disconnected.  Then
by \cite[7.5]{hewittr} there is a neighbourhood basis at $e$ consisting
of compact open subgroups.  Thus for each $y\iin Y$, there is a compact
open subgroup $K_y$ for which $yK_y\subset Y$, so $Y=\bigcup_{y\in Y}yK_y$.
Thus there are $y_1,\dots,y_n\iin Y$ for which $Y=\bigcup_{i=1}^n y_iK_{y_i}$,
so $Y\in\Omega_o^c(G)$.

If $G_e\not=\{e\}$, we have for each $s\in G$ that $sG_e$ is the connected
component of $s$. Hence $Y\supseteq sG_e$ if $s\in Y$ and $Y\cap sG_e=
\varnothing$ otherwise.  Thus $1_Y\in\fal{G\!:\!G_e}$
whence $Y=q_e^{-1}(q_e(Y))$ where $q_e:G\to G/G_e$ is the quotient map, and
$1_{q_e(Y)}\in\fal{G/G_e}$.  Hence $q_e(Y)\in\Omega_o^c(G/G_e)$ from above,
so $Y\in\Omega_o^c(G)$.

{\bf (ii)} We recall that if $G_e$ is compact, then
the map $u\mapsto u\comp q_e:\fal{G/G_e}\to\fal{G:G_e}$
is an isometric isomorphism, so $\fal{G\!:\!G_e}$ is a closed subalgebra of 
$\falg$.  Hence it follows from (i) that $\falig\subseteq\fal{G\!:\!G_e}$.
Conversely, it follows from \cite[Theorem 5.1 and Lemma 5.2]{forrest} that 
$\fal{G\!:\!G_e}\cong\fal{G/G_e}$ is generated by idempotents, so 
$\falig\supseteq\fal{G\!:\!G_e}$.  \endpf

\section{The algebra $\fsalig$}

Let us begin by reviewing some critical ideas and notations from
\cite{ilies2}.  If $G$ is a locally compact group with topology $\tau_G$, we 
let
\[
\tg=\left\{\tau\subseteq\tau_G:
\begin{matrix}\tau\text{ is a group topology on }G\text{ for which} \\
the\text{ completion }G_\tau\text{ is locally compact}
\end{matrix}\right\}
\]
If $\tau\in\tg$, then there is a natural continuous homomorphism
$\eta_\tau:G\to G_\tau$ which is injective if and only if $\tau$
is Hausdorff.  If $\tau_1,\tau_2\in\tg$, and $\tau_1\subseteq\tau_2$,
then there is a continuous homomorphism 
$\eta^{\tau_2}_{\tau_1}:G_{\tau_2}\to G_{\tau_1}$ such that
$\eta^{\tau_2}_{\tau_1}\comp\eta_{\tau_2}=\eta_{\tau_1}$.  Furthermore,
we say that $\tau_1$ is a {\it quotient} of $\tau_2$ if 
$\eta^{\tau_2}_{\tau_1}$ is proper, i.e., for every compact subset $C$ of
$G_{\tau_1}$, $\left(\eta^{\tau_2}_{\tau_1}\right)^{-1}(C)$ is a compact
subset of $G_{\tau_2}$.  This is equivalent to having $\eta^{\tau_2}_{\tau_1}$
be a quotient map with compact kernel.  We then let
\[
\tnqg=\{\tau\in\tg:\tau\text{ is not a proper quotient of some }
\tau_1\iin\tg\}
\]
Every $\tau\iin\tg$ is the quotient of a unique $\tau_{nq}\iin\tnqg$.

For each $\tau\iin\tg$, the map
$u\mapsto u\comp\eta_\tau$ from $\fal{G_\tau}$ to $\fsalg$ is an isometric
homomorphism.  We let
\[
\falbg{\tau}=\fal{G_\tau}\comp\eta_\tau.
\]
This is the space of matrix coefficients of the unitary representation
\begin{equation}\label{eq:lamtau}
\lam_\tau=\lam_{G_\tau}\comp\eta_\tau:G\to\fB(\bltwo{G_\tau})
\end{equation}
where $\lam_{G_\tau}:G_\tau\to\fB(\bltwo{G_\tau})$is the left regular
representation, given by left translations.
We then have that $\falbg{\tau_1}\subseteq\falbg{\tau_2}$ if, and only if,
$\tau_1$ is a quotient of $\tau_2$.  We also have that if $\tau_1$ and $\tau_2$
are distinct elements of $\tnqg$, then
$\falbg{\tau_1}\cap\falbg{\tau_2}=\{0\}$.
Thus we may let the {\it spine} of $\fsalg$ be given by
\[
\spineg=\wbar{\sum_\tau\in\tg}\falbg{\tau}
=\ell^1\text{-}\!\!\!\bigoplus_{\tau\in\tnqg}\falbg{\tau}.
\]
We recall from \cite[Section 2.4]{ilies2} that $\tg$ is a semilattice,
i.e.\ a commutative idempotent semigroup, under the operation
$(\tau_1,\tau_2)\mapsto\tau_1\vee\tau_2$.  The topology $\tau_1\vee\tau_2$
is the coarsest topology on $G$ for which that map 
$s\mapsto\bigl(\eta_{\tau_1}(s),\eta_{\tau_2}(s)\bigr):G\to G_{\tau_1}\cross
G_{\tau_2}$ is continuous.  Thus $G_{\tau_1\vee\tau_2}$ is the closure
of the ``diagonal'' subgroup in $G_{\tau_1}\cross G_{\tau_2}$.
We have that $\tnqg$ is a subsemilattice of $\tg$.
Moreover, $\spineg$ is graded over this semilattice in the sense that
$\falbg{\tau_1}\falbg{\tau_2}\subseteq\falbg{\tau_1\vee\tau_2}$
if $\tau_1,\tau_2\in\tnqg$.

\subsection{The sublattice of pre-totally disconnected topologies}
We say that a locally precompact topology $\tau$ on $G$ is
{\it pre-totally disconnected} if $G_\tau$ is a totally disconnected group,
i.e., its connected component of the identity $G_{\tau,e}$ is trivial.
We let
\[
\ttdg=\{\tau\in\tg:\tau\text{ is pre-totally disconnected}\}.
\]
Then $\ttdg$ is a subsemilattice of $\tg$ for if $\tau_1,\tau_2\in
\ttdg$, then $G_{\tau_1\vee\tau_2,e}$ is a subgroup of
$G_{\tau_1,e}\cross G_{\tau_2,e}$, which is trivial.
It is clear that if $G$ is connected, then
$\ttdg=\{\eps\}$, where $\eps=\{\varnothing, G\}$.  However,
we can say more.

\begin{ttdgisttdge}\label{prop:ttdgisttdge}
There is a natural map $\tau\mapsto\til{\tau}:\ttdg\to
\ttd{G/G_e}$ which is a semilattice isomorphism.
Moreover, $G_\tau=(G/G_e)_{\til{\tau}}$ for each
$\tau\iin\ttdg$.
\end{ttdgisttdge}

\proof If $\tau\in\ttdg$, then the homomorphism
$\eta_\tau:G\to G_\tau$ satisfies $G_e\subseteq\ker\eta_\tau$.
Hence the first isomorphism theorem for topological groups tells 
us that there is a continuous homomorphism
$\til{\eta}_\tau:G/G_e\to G_\tau$ such that $\til{\eta}_\tau\comp q_e
=\eta_\tau$, where $q_e:G\to G/G_e$ is the quotient map.
Let $\til{\tau}$ be the coarsest topology on $G/G_e$ which makes
$\til{\eta}_\tau$ continuous.  Then it is clear that $\til{\tau}\in\ttd{G/G_e}$
and $(G/G_e)_{\til{\tau}}=G_\tau$.  It is straightforward to verify
that $\tau\mapsto
\til{\tau}:\ttdg\to\ttd{G/G_e}$ is bijective and that it is a semilattice
homomorphism.  \endpf

Recall the definition of a quotient topology, given above.
Analogously with the definition of the non-quotient topologies
on $G$, we let
\[
\tnqtdg=\{\tau\in\ttdg:\tau\text{ is not a quotient
of any other element of }\ttdg\}.
\]
We note that the topology induced
by the quotient map $q_e:G\to G/G_e$, is an element of
$\tnqtdg$.  We will identify one other important element
of $\tnqtdg$.  We recall that a {\it compactification} of $G$
is any pair $(\eta,S)$ where $S$ is a compact semi-topological
semigroup and $\eta:G\to S$ is a homomorphism with dense range.
The {\it maximal almost periodic compactification} $(\eta_{ap},G^{ap})$
is the maximal one amongst the (semi-)topological group 
compactifications.  See \cite{berglundjm}, for example.
Note that for a totally disconnected group $G$, even a discrete group,
$G^{ap}$ might not be totally disconnected.  We thus define
the {\it totally disconnected compactification}, $(\eta_{td},G^{td})$,
by
\[
G^{td}=G^{ap}/(G^{ap})_e\aand\eta_{td}= q^{ap}_e\comp\eta_{ap}
\]
where $(G^{ap})_e$ is the connected component of the identity
and $q^{ap}_e:G^{ap}\to G^{ap}/(G^{ap})_e$ is the quotient map.
It follows the first isomorphism theorem and the universal property
for $(\eta_{ap},G^{ap})$
that this compactification has the following universal
property.  

\medskip
\parbox{4.5in}{If $(\eta_1,G_1)$ is a group compactification of $G$
with $G_1$ totally disconnected, then it is a quotient of
$(\eta_{td},G^{td})$ in the sense that there is a quotient
homomorphism $\eta^{td}_1:G^{td}\to G_1$ such that
$\eta^{td}_1\comp\eta_{td}=\eta_1$.}

\begin{tdcompact}\label{prop:tdcompact}
There is an isomorphism $G^{td}\cong(G/G_e)^{td}$.
In fact the compactifications of $G$, $(\eta_{td},G^{td})$
and $\bigl(\til{\eta}_{td}\comp q_e,(G/G_e)^{td}\bigr)$
coincide, where $\bigl(\til{\eta}_{td},(G/G_e)^{td}\bigr)$ is
the totally disconnected compactification of $G/G_e$.
\end{tdcompact}

\proof  Consider the following commutative diagram.
\[
\xymatrix{
G \ar[rr]^{\eta_{td}} \ar[d]_{q_e}
& & G^{td} \ar@{-->}@<1ex>[d]^{\eta_4}  \\
G/G_e \ar@{-->}[urr]^{\eta_1} \ar[r]_{\til{\eta}_{ap}}
& (G/G_e)^{ap} \ar@{-->}[ur]^<<{\eta_2} \ar[r]_{\til{q}_e}
& (G/G_e)^{td} \ar@{-->}@<1ex>[u]^{\eta_3}
}
\]
The homomorphism $\eta_1$ exists by virtue of the
first isomorphism theorem.  The quotient homomorphism
$\eta_2$ exists by the universal property of the compactification
$\bigl(\til{\eta}_{ap},(G/G_e)^{ap}\bigr)$.  Then $\eta_3$ exists,
and is a quotient map, by the first isomorphism theorem.
The inverse, $\eta_4=\eta_3^{-1}$ exists by the universal
property for the totally disconnected compactification
of $G$, applied to the compactification
$\bigl(\til{q}_e\comp\til{\eta}_{ap}\comp q_e,(G/G_e)^{td}\bigr)
=\bigl(\til{\eta}_{td}\comp q_e,(G/G_e)^{td}\bigr)$.
\endpf

We can now give a useful characterisation of the non-quotient
pre-totally disconnected topologies.

\begin{tnqtdgchar}\label{theo:tnqtdgchar}
If $\tau\in\ttdg$, then the following hold.

{\bf (i)} $\tau$ is a quotient of $\tau\vee\tau_{td}$.

{\bf (ii)} If $\tau_1\iin\ttdg$ is any element of which $\tau$
is a quotient, then $\tau_1$ is a quotient of $\tau\vee\tau_{td}$.
Hence $\tau\vee\tau_{td}$ is the unique element of 
$\tnqtdg$ of which $\tau$ is a quotient.

{\bf (iii)} $\tau\in\tnqtdg$ \quad $\iff$ \quad 
$\tau=\tau\vee\tau_{td}$ \quad $\iff$ \quad $\tau\supseteq\tau_{td}$.
\end{tnqtdgchar}

\proof This is proved exactly as \cite[Theorem 2.2]{ilies2}, with $\tau_{td}$
playing the role of the topology induced by the almost periodic
compactification $\tau_{ap}$.   \endpf

\begin{tnqtdgchar1}\label{cor:tnqtdgchar1}
{\bf (i)} $\tnqtdg$ is a subsemilattice of $\ttdg$.

{\bf (ii)} The map $\tau\mapsto\til{\tau}$ of Proposition
\ref{prop:ttdgisttdge} restricts to a semilattice isomorphism
from $\tnqtdg$ to $\tnqtd{G/G_e}$.
\end{tnqtdgchar1}

\proof {\bf (i)} By part (iii) of the theorem above, $\tnqtdg$
is the ideal $\ttdg\vee\tau_{td}$.

{\bf (ii)} It follows Proposition \ref{prop:tdcompact} that
$\wtil{\tau_{td}}=\til{\tau}_{td}$, the topology induced
by the totally disconnected compactification
map $\til{\eta}_{td}:G/G_e\to(G/G_e)^{td}$.  The
result now follows (i) above.  \endpf

\subsection{The structure of $\fsalig$}
We can now state our main structure theorem for $\fsalig$.

\begin{fsalgidem}\label{theo:fsalgidem}
We have that
\[
\fsalig=\ell^1\text{-}\!\!\!\!\bigoplus_{\tau\in\tnqtdg}\faltg
\]
and is a graded Banach algebra over the semilattice $\tnqtdg$,
i.e.\
\[
\falbg{\tau_1}\falbg{\tau_2}\subseteq\falbg{\tau_1\vee\tau_2}
\]
for $\tau_1,\tau_2\iin\tnqtdg$.
\end{fsalgidem}

\proof It was shown in the proof of \cite[Theorem 4.12]{ilies2} that
if $u\in\idemg$, then there are $\tau_1,\dots,\tau_n\iin\tnqg$
and idempotents $u_j\iin\falbg{\tau_j}$ such that
\[
u=u_1+\dots+u_n.
\]
We have by Proposition \ref{prop:falgidem} that
\[
\wbar{\spn}\bigl(\falbg{\tau_j}\cap\idemg\bigr)
\cong\fal{G_{\tau_j}:G_{\tau_j,e}}
\cong\fal{G_{\tau_j}/G_{\tau_j,e}}.
\]
Hence $u_j\in\falbg{\tau'_j}$ where $\tau'_j$ is the topology
induced by the map $q_{\tau_j,e}\comp\eta_{\tau_j}:G\to
G_{\tau_j}/G_{\tau_j,e}$.  Let $\tau''_j=\tau'_j\vee\tau_{td}$.
Then it follows from \cite[Lemma 3.2]{ilies2} that $u_j\in
\falbg{\tau''_j}$.  Hence $u\in\sum_{j=1}^n\falbg{\tau''_j}$.
Thus we have that $\idemg\subset\sum_{\tau\in\tnqtdg}
\faltg$ and hence $\fsalig\subseteq\wbar{\sum_{\tau\in\tnqtdg}
\faltg}$.  Since each $\faltg\cong\fal{G_\tau}$ is generated
by idempotents by Proposition \ref{prop:falgidem}, we obtain
the reverse inclusion, whence
\begin{equation}\label{eq:asum}
\fsalig=\wbar{\sum_{\tau\in\tnqtdg}\faltg}.
\end{equation}

Let us now see that this sum is an $\ell^1$-direct sum.
If $\tau_1\aand\tau_2$ are elements of $\tnqtdg$
and $\falbg{\tau_1}\cap\falbg{\tau_2}\not=\{0\}$, then
by \cite[Lemma 3.2]{ilies2} there is a $\tau\in\tg$
which is a quotient of each of $\tau_1$ and $\tau_2$,
such that
\[
\falbg{\tau_1}\cap\falbg{\tau_2}=\faltg.
\]
It then follows that $\fal{G_\tau}\cong\faltg$ is generated
by idempotents, so $G_\tau$ is totally disconnected, i.e.,
$\tau\in\ttdg$.  It then follows Theorem \ref{theo:tnqtdgchar}
that each of $\tau_1\aand\tau_2$ are quotients of
$\tau\vee\tau_{td}$.  But then $\tau_1=\tau\vee\tau_{td}
=\tau_2$.  Thus if $\tau_1\not=\tau_2$, then $\falbg{\tau_1}\cap
\falbg{\tau_2}=\{0\}$.  It then follows \cite[3.9]{arsac}
that the sum (\ref{eq:asum}) is an $\ell^1$-direct sum.  

The fact that $\fsalig$ is graded over the semilattice
$\tnqtdg$ follows from \cite[Proposition 3.1]{ilies2}.  \endpf

We now aim to describe the Gelfand spectrum of $\fsalig$.
We let
\[
\hdtdg=\left\{\fS\subseteq\tnqtdg:\begin{matrix}
\fS\text{ is hereditary: }\fS\supseteq\fS_{\tau}\text{ for each }\tau\iin
\fS, \aand \\
\fS\text{ is directed: }\tau_1,\tau_2\in\fS\;\implies\;
\tau_1\til{\vee}\tau_2\in\fS\quad
\end{matrix}\right\}
\]
where $\fS_\tau=\{\tau_0\in\tnqtdg:\tau_0\subseteq\tau\}$.  We call sets
$\fS_\tau$ the {\it principal} hereditary directed subsets.  The inverse
mapping system
\[
\{G_\tau,\eta^{\tau_2}_{\tau_1}:G_{\tau_2}\to G_{\tau_1}\,:\,
\tau,\tau_1,\tau_2\in\fS,\tau_1\subseteq\tau_2\}
\]
gives rise to a projective limit group
\[
G_\fS=\underset{\tau\in\fS}{\underleftarrow{\lim}}G_\tau
=\left\{(s_\tau)_{\tau\in\fS}\in \prod_{\tau\in\fS}G_\tau:
\eta^{\tau_2}_{\tau_1}(s_{\tau_2})=s_{\tau_1},\iif \tau_1\subseteq\tau_2\iin
\fS\right\}.
\]
The following description of the Gelfand spectrum of $\fsalig$ 
is similar to the description of that for $\spineg$ in 
\cite[Section 4.1]{ilies2}.  The proof is is the same, so we omit it.

\begin{spectrum}\label{theo:spectrum}
The coproduct $G^I=\bigsqcup_{\fS\in\hdtdg}G_\fS$ is the Gelfand spectrum
of $\fsalig$ via the following identification:  if $s\in G^I$, so
$s=(s_\tau)_{\tau\in\fS}\iin G_\fS$ for some $\fS$, then the character
$\chi_s$ given for $u=\sum_{\tau\in\tnqtdg}u_\tau$ by
\[
\chi_s(u)=\sum_{\tau\in\fS}\hat{u}_\tau(s_\tau)
\]
where $u_\tau\mapsto\hat{u}_\tau:\faltg\to\fal{G_\tau}$ is the
Gelfand transform for each $\tau$.  Moreover, with the weak* topology
and the multiplication
\[
(s_\tau)_{\tau\in\fS_1}(t_\tau)_{\tau\in\fS_2}
=(s_\tau t_\tau)_{\tau\in\fS_1\cap\fS_2}
\]
$G^I$ is a compact semitopological semigroup.
\end{spectrum}

We note that $G^I$ is, in fact, a {\it Clifford semigroup}, i.e.,
a semilattice of groups.  We note that $\hdtdg$ is itself a semilattice
via the operation $(\fS_1,\fS_2)\mapsto\fS_1\cap\fS_2$.
This semilattice is lattice isomorphic to the semilattice
of (central) idempotents in $G^I$, 
\[
\left\{e_\fS=\bigl(\eta_\tau(e)\bigr)_{\tau\in\fS}:\fS\in\hdtdg\right\}.
\]
For every $s=(s_\tau)_{\tau\in\fS}\in G^I$ we let
$s^*=(s_\tau^{-1})_{\tau\in\fS}$, so $s^*s=ss^*=e_\fS$.

It follows from standard von Neumann algebra techniques
(see, for example, \cite[2.2]{arsac}), that the
dual space of $\fsalig$ is the von Neumann algebra
\[
\vnig=\ell^\infty\text{-}\!\!\!\!\bigoplus_{\tau\in\tnqtdg}\vntg.
\]
This von Neumann algebra is generated by the 
representation
\[
\lam_I=\bigoplus_{\tau\in\tnqtdg}\lam_\tau
:G\to\fB\left(\ell^2\text{-}\!\!\!\!\bigoplus_{\tau\in\tnqtdg}\bltwo{G_\tau}
\right)
\]
where $\lam_\tau$ is the representation from (\ref{eq:lamtau}).
Then $G^I$ identifies naturally with the closure of
$\lam_I(G)$ in the weak operator topology (which is why
$G^I$ is a semitopological semigroup).  
We call $(G^I,\lam_I)$ -- or, by abusing notation, just $G^I$ -- 
the {\it idempotent compactification} of $G$.

The set of unitary elements of $G^I$, i.e., those elements
$s$ for which $s^*s=ss^*=e_{\tnqtdg}$, is exactly
$\lam_I(G)$, which, in its weak operator topology, is
bicontinuously isomorphic with $G/G_e$.

\section{Homomorphisms}

\subsection{Completely bounded homomorphisms}
Completely bounded homomorphism between $\falg$ and $\fsalh$
are discussed in detail in \cite{ilies1,ilies2}.  For the purpose of the
discussion below, let us say that a linear map between preduals
of von Neumann algebras, $\Phi:\fM_*\to\fN_*$ is completely bounded
(contractive) if and only if the adjoint map $\Phi^*:\fN\to\fM$ is
completely bounded (contractive).  We leave all other pertinent 
definitions to our other articles \cite{ilies1,ilies2} and the book
\cite{effrosrB}.

We require a definition from \cite[Section 4.4]{ilies2}.  If $G$ and
$H$ are complete topological groups, a {\it continuous piecewise affine}
map from $H$ to $G$ is any map $\alp:Y\to G$, where $Y\in
\Omega_o(H)$, and there exist disjoint subsets $Y_1,\dots,Y_n$
of $Y$, and open cosets $L_1\supseteq Y_1,\dots,L_n\supseteq Y_n$
and affine maps $\alp_j:L_j\to G$ -- i.e., $\alp_j(rs^{-1}t)=
\alp_j(r)\alp_j(s)^{-1}\alp_j(t)$ for all $r,s,t\iin L_j$ -- 
such that $\alp_j|_{Y_j}=\alp|_{Y_j}$, for $j=1,\dots,n$.
If $H$ and $G$ are locally compact groups, then a
continuous piecewise affine map $\alp:Y\subseteq H\to G^I$
is any map for which $\alp(Y)\subseteq G_\fS$ for some $\fS\iin\hdtdg$
and $\alp:Y\to G_\fS$ is continuous piecewise affine.  We say
that $\alp$ is affine, if $Y$ itself is an open coset and $\alp$ an affine map.

\begin{cbhomo}\label{theo:cbhomo}
Let $G$ and $H$ be locally compact groups.
Let $\Phi:\fsalig\to\fsalh$ be a completely bounded homomorphism.
Then there exits a continuous piecewise affine map 
$\alp:Y\subseteq H\to G^I$, which is constant on cosets of $H_e$,
such that
\[
\Phi u(h)=\begin{cases} u\bigl(\alp(h)\bigr) &\text{if }h\in Y, and \\
0 &\text{otherwise}.\end{cases}
\]
Moreover, we have that

{\bf (i)} $\alp$ is affine, if $\Phi$ is completely contractive; and

{\bf (ii)} $\alp(H)\subset\lam_I(G)$, if $\Phi$ is injective.
\end{cbhomo}

\proof  The map $\alp:Y\subseteq H\to G_\fS$, $\fS\in\hdtdg$,
is obtained exactly as in \cite[Theorem 4.9 (ii)]{ilies2}.
We note that $\alp$ is constant on cosets of $H_e$ since $\alp$ is continuous
and $G_\fS$, a subgroup of a product of totally disconnected groups, 
is itself totally disconnected.

Suppose now that $\Phi$ is completely contractive. 
Let $\tau\in\fS$.  Then $\Phi|_{\faltg}:\faltg\to\fsalh$
induces a completely contractive homomorphism
$\Phi_\tau:\fal{G_\tau}\to\fsalh$.  By \cite[Thorem 3.7]{ilies1},
we must have a continuous affine map $\alp_\tau:L_\tau\to
G_\tau$, such that
\[
\Phi_\tau u(h)=\begin{cases}
u\bigl(\alp_\tau(h)\bigr) & \text{if }h\in L_\tau \\
0 &\text{otherwise}. \end{cases}
\]
It follows from the fact that $\fsalig$ is graded over $\tnqtdg$
that for $\tau\subseteq\tau'\iin\fS$, $\eta^{\tau'}_\tau\comp\alp_{\tau'}
=\alp_\tau$, and hence
the open cosets $L_\tau$ must all coincide.  Hence we denote
this open coset by $L$, and we have an affine map $\alp:H\to G_\fS$.
Thus (i) is demonstrated.

Suppose $\Phi$ is injective.  Let us say that there is an open subset
$U$ of $\lam_I(G)$ such that $\alp(Y)\cap U=\varnothing$.
Then we could find $u\in\fal{G:G_e}=\falbg{\til{\tau}_G}$
so that its support $\supp{u}$ is contained in $U$.  But then
it follows that $\Phi u(h)=u\bigl(\alp(h))=0$ for each $h\iin Y$, and
hence $\Phi u=0$.  This contradicts injectivity, whence (ii).
\endpf

We note that if $\fA$ is any Banach algebra which is the closure
of the span of its idempotents, and $\Psi:\fA\to\fsalh$ is any continuous
homomorphism, then $\Psi(\fA)\subseteq\fsalih$.  This applies, in 
particular, to any von Neumann algebra $\fA$, or to $\fA=\fsalig$.
It thus follows from the above theorem that any completely isometric
isomorphism $\Phi:\fsalig\to\fsalih$ induces the map
$\Phi^*|_{H^I}:H^I\to G^I\cup\{0\}$, such that 
\[
\Phi^*\bigl(\lam_I(h)\bigr)=\begin{cases}
\lam_I\bigl(\alp(h)\bigr) &\text{if }h\in Y \\
0 &\text{otherwise.} \end{cases}
\]
Since $\lam_I(H)$ is dense in $H^I$, this determines
$\Phi^*|_{H^I}$, and hence $\Phi$, completely.

\subsection{Isometric isomorphisms}
We do not need to explicitly use operator space structure
if we restrict our consideration to isometric isomorphisms
$\Phi:\fsalig\to\fsalih$.  This is due to the result characterising
linear surjective isometries between C*-algebras due to
Kadison in \cite{kadison}.  This was cleverly adapted by Walter in 
\cite{walter}
to show that the Banach algebra $\fsalg$ is isometrically isomorphic to 
$\fsalh$ if and only if $G\cong H$, as topological groups.
We note that the following result can be immediately deduced from
the theorem above if $G$ is amenable and $\Phi$ is completely
isometric.

\begin{isomorphism}\label{theo:isomorphism}
Let $G$ and $H$ be locally compact groups.
The the following are equivalent:

{\bf (i)} there is a bicontinuous group isomorphism $G/G_e\cong H/H_e$;

{\bf (ii)} there is a bicontinuous semigroup isomorphism $G^I\cong H^I$;

{\bf (iii)} there is an isometric algebra isomorphism $\fsalig\cong\fsalih$.
\end{isomorphism}

\proof {\bf (i) $\implies$ (ii) and (iii).}  It clearly suffices see that
\[
(G/G_e)^I\cong G^I\quad\aand\quad \fsali{G/G_e}\cong\fsalig.
\]
This follows from the fact that there is a natural bijection
$\tau\to\til{\tau}:\tnqtdg\to\tnqtd{G/G_e}$ from part (ii) of Corollary
\ref{cor:tnqtdgchar1}.

{\bf (ii) $\implies$ (i).}   Let $\gam:G^I\to H^I$ be an isomomorphism.
Then $\gam$ must take the unit, $E_G=\lam_I(e_G)$, of $G^I$
 to that, $E_H=\lam_I(e_H)$, of $H^I$.  Indeed, $\gam$ is surjective
 and $\gam(E_G)$ acts as the unit on $H^I=\gam(G^I)$.
 Now if $U$ is a unitary
 in $G^I$, i.e., $U=\lam_I(s)$ for some $s\iin G$, then
 \[
 \gam(U^*)\gam(U)=\gam(U^*U)=\gam(E_G)=E_H.
 \]
Then $\gam(U)$ is a unitary in $H^I$,
since the idempotent $E_U=\gam(U)^*\gam(U)$ in $H^I$, 
satisfies $E_H=\gam(U^*)\gam(U)=\gam(U^*)\gam(U)E_U=
E_HE_U=E_U$.
Thus $\gam(U)=\lam_I(h)$ for some $h\iin H$, and 
$\gam(U^*)=\gam(U)^*=\lam_I(h^{-1})$.  

Thus, the isomorphism
$\gam:G^I\to H^I$ must satisfy that $\gam\bigl(\lam_I(G)\bigr)
=\lam_I(H)$, and hence must induce an isomorphism
$G/G_e\cong H/H_e$.  Moreover, this isomorphism is bicontinuous
if $\gam$ is.

{\bf (iii) $\implies$ (i).}  
If $\Phi:\fsalig\to\fsalih$ is an isometric isomorphism,
then its adjoint $\Phi^*:\vnih\to\vnig$ is a linearly isomorphic isometry.
Moreover, $\Phi^*(H^I)=G^I$, by Theorem \ref{theo:spectrum}.  
It follows 
\cite[Theorem 1]{kadison} that the unitary subgroup
$\lam_I(H)$ of $H^I$, which is bicontinuously isomorphic to $H/H_e$,
is exactly the set of extreme points of the unit ball of $\vnih$ which lie in
$H^I$.  The same is true for the
unitary subgroup $\lam_I(G)$ of $G^I$.  Hence $\Phi^*\bigl(\lam_I(H)\bigr)
=\lam_I(G)$ and thus $\Phi^*$ induces a homeomorphism $\vphi:H/H_e\to G/G_e$.
We note that $\lam_I(H)$ is a linearly independent subset 
of $\vnih$, as it represents the Gelfand spectrum of the regular subalgebra
$\fal{H:H_e}\cong\fal{H/H_e}$ of $\fsalih$.
It can now be argued, exactly as in the proof of \cite[Theorem 2]{walter}
that $h\mapsto\vphi\bigl(h\vphi(e_H)^{-1}\bigr)$ 
is either an isomorphism or an anti-isomorphism
between the groups $H/H_e$ and $G/G_e$.
(It will be an isomorphism if $\Phi$ is completely contractive.)
Composing this map with the canonical continuous anti-isomorphism
$s\mapsto s^{-1}$ on $G/G_e$, if necessary, we see that
$H/H_e\cong G/G_e$ as topological groups.  \endpf

We cannot weaken the assumptions of part (iii) of the above theorem
too much.  If $G$ and $H$ are non isomorphic finite groups of the
same cardinality then $\fsalig=\falg\cong\fal{H}=\fsalih$ via a completely
bounded isomorphism. 

\subsection{Amenability} We make brief note of some of the amenability
properties of $\fsalig$.  For sake of brevity, and since we will make no
use of these concepts outside of this section,
we refer the reader to \cite{rundes} for all 
pertinent definitions on {\it amenability}, {\it weak amenability}
and {\it operator amenability} of a {\it completely contractive Banach
algebra}.

First, we note that it is immediate from remarks following 
\cite[Theorem 1.4]{badecd} that $\fsalig$, being an algebra generated by 
idempotents, is always weakly amenable, and hence always operator weakly
amenable.  The following result can be proved exactly as
\cite[Theorem 3.3]{rundes} and the corollary which follows it.

\begin{amenability}
Let $G$ be a locally compact group.

{\bf (i)} The completely contractive Banach algebra $\fsalig$
is operator amenable if and only if $G/G_e$ is amenable and
$\tnqtdg$ is finite.

{\bf (ii)} The Banach algebra $\fsalig$ is amenable if and only if
$G/G_e$ has an abelian subgroup of finite index, and
$\tnqtdg$ is finite.
\end{amenability}

\section{Examples}

\subsection{Almost connected groups}  
If $G$ is an almost connected group, i.e., $G/G_e$ is compact,
then there is only one non-quotient topology amongst the
pre-totally disconnected locally pre-compact topologies.
It is the topology induced by the quotient map $q_e:G\to G/G_e$.
Thus we have that
\[
\fsalig=\fsal{G:G_e}\cong\fal{G/G_e}
\]
where $\fsal{G:G_e}$ is the subalgebra of $\fsalg$ of elements
which are constant on cosets of $G_e$.

\subsection{Some Fell groups}  We now give an example of a family
of non-compact groups $G$ for which $\fsalig=\fsalg$.
Let $\Que_p$ denote the field of $p$-adic numbers ($p$ prime) and
$\Oh_p$ the compact subring of $p$-adic integers.  We let
$\genlin{n,\Oh_p}$ denote the compact group of invertible
$n\cross n$ matricies with entries in $\Oh_p$ and each having inverse
with entries in in $\Oh_p$.  These groups were introduced in
\cite{maucerip} for general $n$; and in \cite{bagget} for $n=1$.
We let $\genlin{n,\Oh_p}$ act on the the elements of the
vector space $\Que_p$ as column vectors.  Then we let
\[
G_{p,n}=\genlin{n,\Oh_p}\ltimes\Que_p^n.
\]
It was noted in \cite{rundes} that for each prime $p$ and
$n=1,2,\dots$ that
\[
\fsal{G_{p,n}}=\fsal{G_{p,n}:\Que_p^n}+\fal{G_{p,n}}
\]
where $\fsal{G_{p,n}:\Que_p^n}\cong\fal{G_{p,n}/\Que_p^n}
\cong\fal{\genlin{n,\Oh_p}}$.  Since all groups involved
are totally disconnected, we obtain 
\[
\fsali{G_{p,n}}=\fsal{G_{p,n}}.
\]
This is the only family of non-compact groups for which this is known.

\subsection{The integers $\Zee$ and the rationals $\Que$}
For any locally compact abelian group $G$ we let $\what{G}$ denote
its dual group.
Using the well-known formula $\Zee^{ap}\cong\what{\Tee_d}$, where $\Tee_d$
is the circle group with discrete topology, and then applying
\cite[25.1]{hewittr}, we obtain that $\Zee^{ap}\cong\what{\Que}^\mathfrak{c}
\cross\prod_{p\text{ prime}}\Oh_p$. We let $\mathfrak{c}$ denote the 
cardinality of the real numbers, and as above, $\Oh_p$ the $p$-adic
integers, which is the closure of $\Zee$ in the $p$-adic number field $\Que_p$.
Since $\what{\Que}$ is connected -- it contains a dense copy of the real line,
we obtain the totally disconnected compactification
\[
(\eta_{td},\Zee^{td})=\left(n\mapsto(n,n,\dots),\prod_{p\text{ prime}}\Oh_p
\right).
\]
It then follows from \cite[9.1]{hewittr} that $\tnqtd{\Zee}=\{\tau_d,
\tau_{td}\}$ where $\tau_d$ is the discrete topology.  Hence
\[
\fsali{\Zee}=\falb{\tau_{td}}{\Zee}\oplus_{\ell_1}
\fal{\Zee}.
\]
Each of the idempotents $1_{n\Zee}$, $n=1,2,\dots$, lies in 
$\falb{\tau_{td}}{\Zee}$.

We note that by \cite[25.2 and 23.21]{hewittr}
$(\prod_{p\text{ prime}}\Oh_p)^\wedge\cong
\bigoplus_{p\text{ prime}}\mathrm{t}_p(\Tee)$, where 
$\mathrm{t}_p(\Tee)$ is the $p$-torsion subgroup of $\Tee$, the sum of
which is the torsion subgroup $\mathrm{t}(\Tee)$.  
Applying the Fourier-Stieltjes transform $\fsal{\Zee}\cong\mathrm{M}(\Tee)$
we realize the closed span of the idempotent measures on $\Tee$ as
\[
\mathrm{M}_I(\Tee)=\ell^1(\mathrm{t}(\Tee))\oplus_{\ell_1}
\mathrm{L}^1(\Tee).
\]

We have that $\what{\Que}$ is a divisible torsion-free group by 
\cite[24.23]{hewittr}.  Thus it is a vector space over $\Que$ which,
being connected, has cardinality $\mathfrak{c}$.  Hence the discrete group
$\what{\Que}_d$ is isomorphic to $\Que^{\oplus\mathfrak{c}}$, and 
so we obtain that $\Que^{ap}\cong\what{\what{\Que}_d}\cong
\what{\Que}^\mathfrak{c}$.  Thus $\Que^{ap}$ is connected from which it follows
that $\Que^{td}$ is a one-element group.  It thus follows that
$\ttd{\Que}=\tnqtd{\Que}$.  $\ttd{\Que}$ contains each of the topologies
$\tau_p$, induced from the embedding $\Que\hookrightarrow\Que_p$.
The weak approximation theorem from number theory (see, for example
\cite[Section 6]{cassels}) tells us that the topologies
$\{\tau_p:p\text{ prime}\}$ form a free semilattice within $\ttd{\Que}$.
There are other topologies.  For example, \cite[Lemma 4.10]{ilies2}
shows how to construct a topology $\tau$ on $\Que$ for which
the closure of $\Zee$ in $\Que_\tau$ is an open subgroup isomorphic
to $\Zee^{td}$; this is related to the construction of the
adelic ring (see \cite[Section 14]{cassels}) of $\Que$.
This construction can be repeated with any subgroup of $\Que$.

\subsection{The lattice $\Zee^2$}
For any locally compact group $G$ and $\tau\iin\tnqtdg$, there
is a $\tau_{nq}\iin\tnqg$ of which $\tau$ is a quotient.  Moreover,
$\tau_{nq}$ must satisfy that $G_{\tau_{nq}}$ has compact connected
component of the identity.

The semilattice $\tnq{\Zee^2}$ is computed explicitly in 
\cite[Section 6.1]{ilies2}.  It is isomorphic to the subspace
lattice of $\Ree^2$.  Let us recall this identification, but with
different notation from what is used in \cite{ilies2}.
If $0\leq\theta<2\pi$ let $\tau_\theta$ be the topology
on $\Zee^2$ induced by the map 
\[
\Zee^2\to\Ree:(n,m)\mapsto n\cos\theta+m\sin\theta.
\]
If $\tan\theta$ is irrational $\Zee^2_{\tau_\theta}\cong\Ree$,
otherwise it is isomorphic to a discrete subgroup of $\Ree$.
Then $\tnq{\Zee^2}=\{\tau_{ap},\tau_\theta\vee\tau_{ap},\tau_d:
0\leq\theta<2\pi\}$, where $\tau_d$ is the discrete topology and
$\tau_{ap}$ the maximal pre-compact topology.  
Consider for each $\tau$ in $\tnq{\Zee^2}$ the topology
induced by the map $q_{\tau,e}\comp\eta_\tau:\Zee^2\to
\Zee^2_\tau/\Zee^2_{\tau,e}$.  Collecting the distinct such topologies
we obtain
\[
\tnqtd{\Zee^2}=\bigl\{\tau_{td},\tau_\theta\vee\tau_{td},\tau_d:
0\leq\theta<2\pi\aand\tan\theta\in\{\Que,\infty\}\bigr\}.
\]
The algebras $\mathrm{A}_{\tau_\theta\vee\tau_{td}}(\Zee^2)$, 
contain all of the idempotents
of the form $1_H$, where $H$ is an infinite subgroup of infinite index
in $\Zee^2$.  The algebra $\mathrm{A}_{\tau_{td}}(\Zee^2)$ contains
all of the idempotents of the form $1_L$, where $L$ is a subgroup of finite
index in $\Zee^2$.

\subsection{Free groups}
Let $\mathrm{F}_\infty$ be the countable free group on infinitely many free 
generators $\{x_n:n\in\En\}$.  If $G$ is any separable totally disconnected
and $S=\{s_n:n\in\En\}$ is any subset which generates a dense subgroup,
then the homomorphism
\[
\eta_S:\mathrm{F}_\infty\to G:x_n\mapsto s_n
\]
induces an element of $\ttd{\mathrm{F}_\infty}$, namely
$\tau_S=\{\eta_S^{-1}(U):U\text{ is open in }G\}$.
Similarly, $\mathrm{F}_\infty^{td}$ 
admits every totally disconnected separable group as a quotient.
We expect $\tnqtd{\mathrm{F}_\infty}$ to be very large.  If we let
$\mathrm{F}_n$ ($n\geq 2$) be a finitely generated free group, then
it admits a copy of  $\mathrm{F}_\infty$ as an (open) subgroup.
Using an adaptation of \cite[Lemma 4.10]{ilies2} we can, for every
$\tau\iin\tnqtd{\mathrm{F}_\infty}$, construct a topology
$\bar{\tau}$ in $\tnqtd{\mathrm{F}_n}$ for which there is an open subgroup 
which is a copy of $(\mathrm{F}_\infty)_\tau$.

{
\bibliography{idempotentbib}

\begin{thebibliography}{10}

\bibitem{arsac}
G.~Arsac.
\newblock Sur l'espace de {B}anach engendr\'{e} par les coefficients d'une
  repr\'{e}sentation unitaire.
\newblock {\em Pub. D\'{e}p. Math. Lyon}, 13(2):1--101, 1976.

\bibitem{badecd}
W.~G. Bade, P.~C. Curtis, and H.~G. Dales.
\newblock Amenability and weak amenability for {B}eurling and {L}ipschitz
  algebras.
\newblock {\em Proc. London Math. Soc.}, 55:359--377, 1987.

\bibitem{bagget}
L.~Bagget.
\newblock A separable group having discrete dual space is compact.
\newblock {\em J. Funct. Anal.}, 10:131--148, 1972.

\bibitem{berglundjm}
J.~F. Berglund, H.~D. Junghelm, and P.~Milnes.
\newblock {\em Analysis on Semigroups: Function Spaces, Compactifications,
  Representations}.
\newblock Canadian Math.\ Society Monographs and Advanced Texts.
  Wiley-Interscience, New York, 1989.

\bibitem{cassels}
J.~W.~S. Cassels.
\newblock Global fields.
\newblock In {\em Algebraic number theory. Proceedings of an instructional
  conference organized by the London Mathematical Society}, pages 42--84,
  London, 1967. Academic Press.

\bibitem{effrosrB}
E.~G. Effros and Z.-J. Ruan.
\newblock {\em Operator Spaces}, volume~23 of {\em London Math. Soc., New
  Series}.
\newblock Claredon Press, Oxford Univ. Press, New York, 2000.

\bibitem{eymard}
P.~Eymard.
\newblock L'alg\`{e}bre de {F}ourier d'un groupe localement compact.
\newblock {\em Bull. Soc. Math. France}, 92:181--236, 1964.

\bibitem{forrest}
B.E. Forrest.
\newblock Fourier analysis on coset spaces.
\newblock {\em Rocky Mountain J. Math.}, 28:173--190, 1998.

\bibitem{hewittr}
E.~Hewitt and K.~A. Ross.
\newblock {\em Abstract Harmonic Analysis I}, volume 115 of {\em Grundlehern
  der mathemarischen Wissenschaften}.
\newblock Springer, New York, second edition, 1979.

\bibitem{host}
B.~Host.
\newblock La th\'{e}or\`{e}me des idempotents dans ${B}({G})$.
\newblock {\em Bull. Soc. Math. France}, 114:215--223, 1986.

\bibitem{ilies1}
M.~Ilie and N.~Spronk.
\newblock Completely bounded homomorphisms of the {F}ourier algebras.
\newblock {\em J. Funct. Anal.}, 2005.
\newblock Accepted.

\bibitem{ilies2}
M.~Ilie and N.~Spronk.
\newblock The spine of a {F}ourier-{S}tieltjes algebra.
\newblock Preprint, see Arxiv {\tt math.FA/0405063}, 2005.

\bibitem{kadison}
R.~V. Kadison.
\newblock Isometries of operator algebras.
\newblock {\em Ann. Math.}, 54:325--338, 1951.

\bibitem{maucerip}
G.~Mauceri and M.~A. Picardello.
\newblock Non-compact unimodular groups with purely atomic {P}lancherel
  measures.
\newblock {\em Proc. Amer. Math. Soc.}, 78:77--84, 1980.

\bibitem{rundes}
V.~Runde and N.~Spronk.
\newblock Operator amenability of {F}ourier-{S}tieltjes algebras, {II}.
\newblock Preprint, see ArXiv {\tt math.FA/0507373}, 2005.

\bibitem{taylor}
J.L. Taylor.
\newblock {\em Measure Algebras}, volume~16 of {\em Conference Board of
  Mathematical Sciences}.
\newblock American Mathematical Society, Providence, RI, 1973.

\bibitem{walter}
M.~E. Walter.
\newblock W*-algebras and nonabelian harmonic analysis.
\newblock {\em J. Func. Anal.}, 11:17--38, 1972.

\end{thebibliography}
\bibliographystyle{plain}
}

\smallskip
{\sc Department of Mathematical Sciences, Lakehead University, 955
Oliver Road, Thunder bay, ON\quad P7B 5E1, Canada}

E-mail: {\tt milie@lakeheadu.ca}

\smallskip
{\sc Department of Pure Mathematics, University of Waterloo,
Waterloo, ON\quad N2L 3G1, Canada} 

E-mail: {\tt nspronk@uwaterloo.ca}

\end{document}